\documentclass[12pt]{article}
\usepackage{amscd,amssymb,amsmath,amsthm}

\newcommand\cD{{\cal{D}}}

\newcommand\cF{{\cal{F}}}
\newcommand\cG{{\cal{G}}}

\newcommand\cL{{\cal{L}}}

\newcommand\cT{{\cal{T}}}

\newcommand\cV{{\cal{V}}}

\newcommand\cZ{{\cal{Z}}}

\newcommand\RR{{\mathbb{R}}}

\newcommand\NN{{\mathbb{N}}}

\theoremstyle{plain}
\newtheorem{thm}{Theorem}
\newtheorem{proposition}[thm]{Proposition}
\newtheorem{lemma}[thm]{Lemma}

\newtheorem{cor}[thm]{Corollary}
\newtheorem{lem}[thm]{Lemma}

\newtheorem{prop}[thm]{Proposition}

\theoremstyle{definition}
\newtheorem{defn}[thm]{Definition}
\newtheorem{ex}[thm]{Example}

\newtheorem{main:as}{Main Assumption}

\theoremstyle{remark}

\numberwithin{equation}{section}

\begin{document}
\title{The trace formula for transversally elliptic
operators on Riemannian foliations}
\author{Yuri A. KORDYUKOV\\
\footnotesize{Department of Mathematics,}\\
\footnotesize{Ufa State Aviation Technical University,}\\
\footnotesize{K.Marx str. 12,
Ufa 450000, Russia}\\ \footnotesize{e-mail: yurikor@math.ugatu.ac.ru}}
\date{September 1, 1999}
\maketitle



\bibliographystyle{unsrt}

\section{Introduction}
The main goal of the paper is to generalize the
Duistermaat-Guillemin trace formula to the case of
transversally elliptic operators on a compact foliated manifold.
First, let us recall briefly the setting of the classical formula.

Let $P$ be a positive self-adjoint elliptic pseudodifferential operator
of order one
on a closed manifold $M$ (for example, $P=\sqrt{\Delta}$, where $\Delta$ is the
Laplace-Beltrami operator of a Riemannian metric on $M$).
For any function $f\in C^{\infty}_c(\RR)$, the operator
$U_f=\int f(t)e^{itP}dt$ can be shown to be of trace class, and
the mapping $\theta:f\mapsto {\rm tr}\,U_f$ is a continuous linear
functional on $C^{\infty}_c(\RR)$. Otherwise speaking,
$\theta$ is a distribution on $\RR$.
The principal symbol $p$ of the operator $P$
is a smooth function on the symplectic manifold $T^*M\setminus 0$,
and, by definition, the bicharacteristic flow $f_t$ of the
operator $P$ is the Hamiltonian flow on $T^{\ast}M\setminus 0$
defined by the function $p$.

By the theorem due to Colin de Verdi\'{e}re and Chazarain \cite{CdV73,Chaz},
the singularities of the distribution $\theta$ are contained in the period set
of closed trajectories of the bicharacteristic flow $f_t$.
Moreover, Duistermaat and Guillemin showed \cite{DG} that,
under the assumption that the bicharacteristic flow
is clean, one can write down an asymptotic expansion for the
distribution $\theta$ near a given period of closed bicharacteristic.
A formula for the leading term of this asymptotic expansion
is the Duistermaat-Guillemin trace formula mentioned above.
It involves the geometry of the bicharacteristic flow in
the form of Poincar\'e map and Maslov indices and
provides a far-reaching generalization of the classical Poisson formula
and the Selberg trace formula on hyperbolic spaces.

In the paper, we prove a trace formula for an operator $P=\sqrt{A}$, where
$A$ is a positive self-adjoint transversally elliptic
pseudodifferential operator of order two with the positive, holonomy invariant
transverse principal symbol on a compact foliated manifold $(M,\cF)$
(see Theorem~\ref{t:main}).
One can consider such an operator as an elliptic operator on
the singular space $M/\cF$ of leaves of the foliation $\cF$
(this statement can be made more
precise, using the language of noncommutative geometry, see \cite{noncom})
the trace formula stated in this paper as an example of a
trace formula for elliptic operators on singular spaces. We hope that
this formula will be useful in further study of a general trace formula
in noncommutative geometry
(see, for instance, \cite{Connes:zeta,Golse:Leichtnam} for discussion of a
noncommutative trace formula).

It should be also noted that our trace formula can be viewed as a relative
version of the Duistermaat-Guillemin trace formula.

\section{Preliminaries and main results}
Let $(M,{\cal F})$ be a compact, connected, oriented foliated manifold.
We will use the following notation: $T\cF$ is the tangent bundle,
$H\cF=TM/T\cF$ is the normal bundle and $N^*{\cal F}$ is the conormal bundle
to ${\cal F}$.
There is a short exact sequence
\begin{equation}\label{e:short}
0\longrightarrow T\cF \longrightarrow TM \longrightarrow H\cF
\longrightarrow 0.
\end{equation}

We will consider linear operators, acting on half-densities. Recall that an
$\alpha$-density ($\alpha\in \RR$) on a real vector space $V$ of dimension
$n$ is a map
$\phi:\Lambda^nV\rightarrow \RR$ such that $\phi(\lambda v)=
|\lambda|^{\alpha}\phi(v), v\in \Lambda^nV, \lambda\in\RR$.
For any real vector bundle $E$ over a smooth manifold $X$, we will denote
by $|E|^{\alpha}$ the $\alpha$-density bundle of $E$.

Given a pseudodifferential operator $A\in  \Psi ^{m}(M,|TM|^{1/2})$,
{\bf the transversal principal symbol} $\sigma_{A}$ of $A$
is defined to be the restriction of its principal symbol $a_m$
on $\tilde{N}^{\ast}{\cal F}=N^*{\cF} \setminus 0$.
An operator $A\in\Psi^m(M,|TM|^{1/2})$ is said to be {\bf transversally
elliptic}, if
$\sigma_A(\nu)\not=0$ for any $\nu\in \tilde{N}^{\ast}{\cal F}$.

For any smooth leafwise path $\gamma$ from $x\in M$ to $y\in M$,
sliding along leaves of the foliation defines {\bf the holonomy map}
$h_{\gamma}$,
which associates to every germ of a local transversal to the foliation
at the point $x$
a germ of a local transversal to the foliation at the point $y$ (this map is
a natural generalization
of the Poincar\'e first-return map for flows). The differential of this map
(the linear holonomy map)
is well-defined as a linear map $dh_{\gamma}:H_x{\cF}\to H_y{\cF}$
and the codifferential as a linear map $dh_{\gamma}^{*}:N^{*}_{y}{\cal F}
\rightarrow
N^{*}_{x}{\cal F}$.

The transversal principal symbol $\sigma_A$ of an operator
$A\in\Psi^m(M,|TM|^{1/2})$
is said to be {\bf holonomy invariant}, if
$\sigma _{A}(dh_{\gamma}^{\ast}(\nu))
= \sigma _{A}(\nu )$ for any  smooth  leafwise path $\gamma $ from $x$ to $y$
and for any $\nu  \in \tilde{N}^{ \ast}_{y}{\cal F}$.
\medskip\par
Throughout in the paper, we will assume that $A$ is a linear operator in
$C^{\infty}(M,|TM|^{1/2})$,  satisfying the following conditions:
\bigskip\par
(A1) $A\in \Psi^2(M,|TM|^{1/2})$ is a transversally elliptic operator with
the positive, holonomy invariant transversal principal symbol;
\medskip\par
(A2) $A$ is an essentially self-adjoint positive operator in $L^2$ space of
half-densities on $M$, $L^2(M)$ (with the
initial domain $C^{\infty}(M,|TM|^{1/2})$).

\begin{ex}
A geometrical example of an operator, satisfying the conditions (A1) and (A2),
is given by the operator $A=I+\Delta_H$, where $\Delta_H$ is the transversal
Laplacian
of a bundle-like metric on a Riemannian foliation.

Recall that a foliation $\cF$ on a smooth Riemannian manifold $(M,g_M)$
is {\bf Riemannian} if it satisfies one of the following equivalent
conditions (see, for instance, \cite{Re}):
\begin{enumerate}
\item $(M,{\cal F})$ locally has the structure of Riemannian submersion;
\item the transverse part of the Riemannian metric $g_M$ (that is,
its restriction to $H=T\cF^{\bot}$) is holonomy invariant;
\item the horizontal distribution $H$ is totally geodesic.
\end{enumerate}
In this case, the metric $g_M$ is called {\bf bundle-like}.

The Riemannian metric $g_M$ defines a decomposition of the
cotangent bundle $T^*M$ into a direct sum
$T^*M=F^*\oplus H^*$. With respect to this decomposition,
the de Rham differential $d:C^{\infty}(M)\to C^{\infty}(M,T^*M)$ can be written
as a sum $d=d_F+d_H,$ where $d_F:C^{\infty}(M)\to C^{\infty}(M,F^*)$ and
$d_H:C^{\infty}(M)\to C^{\infty}(M,H^*)$.

The transversal Laplacian is a second order transversally elliptic
differential operator
in the space $C^{\infty}(M)$ defined by the formula
$$
\Delta_H=-d^*_Hd_H.
$$
Its principal symbol $a_2$ is given by the formula
$$
a_2(x,\xi)=g_H(\xi,\xi), \quad (x,\xi)\in \tilde{T}^*M,
$$
and the holonomy invariance of the transverse
principal symbol $\sigma_{\Delta_H}$ is equivalent
to the assumption on the Riemannian metric $g_M$ to be bundle-like.
\end{ex}

From now on, we will assume that $A$ satisfies the assumptions (A1) and (A2).
By the spectral theorem, the operator $P=\sqrt{A}$ generates
a strongly continuous group $e^{itP}$ of bounded operators in $L^2(M)$.
To define a distributional trace of the operator $e^{itP}$, one need
an additional regularization. First, let us introduce some notation.

Recall that the holonomy groupoid $G=G_{\cF}$ of the foliation ${\cal F}$
is the set of equivalence classes
of leafwise paths $\gamma:[0,1]\rightarrow M$ with respect to
an equivalence relation $\sim_h$, setting $\gamma_1\sim_h \gamma_2$
if $\gamma_1$ and $\gamma_2$ have the same initial
and final points and the same holonomy maps.
$G$ is equipped with maps $s,r:G\rightarrow M$ given by $s(\gamma)=\gamma(0)$
and
$r(\gamma)=\gamma(1)$ and has a composition law given by the composition
of paths. For any $\gamma_1,\gamma_2\in G$, the composition $\gamma_1\circ
\gamma_2$ makes sense iff $r(\gamma_2)=s(\gamma_1)$.
We will make use of standard notation:
$G^x=r^{-1}(x)$, $G_x=s^{-1}(x)$,
$G^x_x=s^{-1}(x)\cap r^{-1}(x)$, $x\in M$.
For any $x\in M$, $G^x_x$ is the holonomy group of the leaf $L_x$
through the point $x$ and the maps $s:G^x\rightarrow L_x$ and
$r:G_x\rightarrow L_x$ are covering maps associated with $G^x_x$.
We will identify a point $x\in M$ with the element in $G$ given
by the constant path $\gamma(t)=x, t\in [0,1]$.

Let $s^*(|T{\cal F}|^{1/2})$ and $r^*(|T{\cal F}|^{1/2})$
be the lifts of the vector bundle of leafwise half-densities
$|T{\cal F}|^{1/2}$
to vector bundles on $G$ via the mappings $s$ and
$r$ respectively, and $|T{\cG}|^{1/2}=r^*(|T{\cal F}|^{1/2})\otimes
s^*(|T{\cal F}|^{1/2}).$ The line bundle $|T{\cG}|^{1/2}$ is the bundle
of leafwise half-densities on $G$ with respect to the natural foliation
$\cG$ \cite{Co79}.

The space $C^{\infty}_c(G,|T{\cG}|^{1/2})$ has the structure of involutive
algebra
(see, for instance, \cite{Co79}).
There is a natural $\ast$-representation $R$
of the involutive algebra $C^{\infty}_c(G,|T{\cG}|^{1/2})$ in $L^2(M)$.
For any $k\in C^{\infty}_c(G,|T{\cG}|^{1/2})$, the operator $R(k)$ in
$L^2(M)$ is defined as follows.
According to the short exact sequence \eqref{e:short},
the half-density vector bundle $|TM|^{1/2}$ can be decomposed as
$$
|TM|^{1/2}\cong |T{{\cal F}}|^{1/2} \otimes |H{\cal F}|^{1/2}.
$$
For any $\gamma\in G, s(\gamma)=x, r(\gamma)=y$,
the corresponding linear holonomy map defines a map
$$dh^*_{\gamma}: |H_y{\cal F}|^{1/2}\to |H_x{\cal F}|^{1/2}.$$
Given $u\in L^2(M)$ of the form $u=u_1\otimes u_2, u_1\in
L^{2}(M,|T\cF|^{1/2}),
u_2\in L^{2}(M,|H\cF|^{1/2}),$
$R(k)u\in L^{2}(M)$ is defined by the formula
$$
R(k)u(x) =\int_{ G^{x}}  k(\gamma )s^*u_1(\gamma )\otimes
dh^*_{\gamma^{-1}}[u_2(s(\gamma))],\quad  x\in  M.
$$

\begin{prop}\label{p:existence}
For any $k\in C^{\infty}_c(G,|T{\cG}|^{1/2})$ and $f\in C^{\infty}_c({\RR})$,
the operator $R(k)\int f(t)e^{itP}dt$ is of trace class.
Moreover, for any $k\in C^{\infty}_c(G,|T{\cG}|^{1/2})$, the formula
$$
\langle\theta_k,f\rangle= {\rm tr}\, R(k)\int f(t)e^{itP}dt, \quad
f\in C^{\infty}_c({\RR}),
$$
defines a distribution $\theta_k$ on the real line $\RR$,
$\theta_k\in {\cal D}'({\RR})$.
\end{prop}

Let ${\cal F}_N$ be a foliation in $\tilde{N}^*{\cal F}$, which is
the horizontal foliation for
the natural leafwise flat connection in $\tilde{N}^*{\cal F}$
(the Bott connection).
The leaf of the foliation $\cF_N$ through a point
$\nu\in \tilde{N}^*{\cal F}$ is the set of all
$dh_{\gamma}^{*}(\nu)\in \tilde{N}^*{\cal F}$
such that $\gamma\in G, r(\gamma)=\pi(\nu)$,
where $\pi:N^*{\cF}\rightarrow M$ is
the natural projection.

Denote by $H{\cal F}_N$ the normal bundle to ${\cF}_N$,
$H{\cal F}_N=T(N^*{\cF})/T{\cF}_N$.
For any $\nu\in N^*{\cal F}$, the space $T_{\nu}N^*{\cal F}$ is a coisotropic
subspace of
the space $T_{\nu}T^*M$ equipped with the canonical symplectic structure,
and $T_{\nu}{{\cal F}_N}$ its skew-orthogonal complement, therefore,
the normal bundle $H_{\nu}{{\cal F}_N}$
has a natural symplectic structure (see, for instance, \cite{GS79}).

Given an operator $A$ under the conditions (A1) and (A2) with the principal
symbol
$a$, let ${\tilde p}$ be a smooth function on $\tilde{T}^*M$
homogeneous of degree one such that $\tilde{p}(\xi)\not=0$ for
$\xi\in \tilde{T}^*M$, which is equal to $p=a^{1/2}$ in some conical
neighborhood
of $N^*{\cal F}$, and $\tilde{f}_t$ the Hamiltonian flow of the
function $\tilde{p}$. Define $\sigma_P$ to be the restriction of
$p$ on $N^*{\cF}$: $\sigma_P=\sigma_A^{1/2}$.
The function $\sigma_P$ coincides with the transverse principal symbol
of any operator $P_1\in \Psi^1(M,|TM|^{1/2})$
such that the principal symbols of $P^2_1$ and $A$ are equal on $N^*{\cal F}$.

The holonomy invariance assumption on $\sigma_A$ implies
\begin{equation}\label{e:hol-inv}
d\tilde{p}(\nu)(X)=0,\quad \nu\in \tilde{N}^*{\cF},\quad X\in T_{\nu}\cF_N.
\end{equation}
Using \eqref{e:hol-inv} and the fact that $\tilde{f}_t$ preserves the
symplectic
structure of $T^*M$, one can easily check that the Hamiltonian flow
$\tilde{f}_t$ can be restricted on $N^*\cF$. The resulting flow
will be denoted by $f_t$. By definition, the flow $f_t$ depends only on
the $1$-jet of the principal symbol $a$ on $N^*{\cF}$,
therefore, it doesn't depend on a choice of $\tilde{p}$
and can be naturally called {\bf the transverse bicharacteristic flow} of
the operator $A$.

Since $\tilde{f}_t$ preserves the symplectic structure of $T^*M$ and
$T{\cF}_N$ is the skew-adjoint complement to $TN^*{\cF}$,
$f_t$ maps leaves of the foliation ${\cal F}_N$
to leaves. In particular, the differential $df_t$ defines
a map $T_{\nu}{{\cal F}_N}\rightarrow T_{f_t(\nu)}{{\cal F}_N}$ and
a symplectic map $H_{\nu}{{\cal F}_N}
\rightarrow H_{f_t(\nu)}{{\cal F}_N}$.

We say that a point $\nu\in \tilde{N}^*{\cal F}$ is {\bf a relative
fixed point} of the
diffeomorphism $f_t:\tilde{N}^*{\cal F}\to \tilde{N}^*{\cal F}$
(with respect to the foliation $\cF_N$), if there exist
$\gamma\in G$ such that $r(\gamma)=\pi(\nu)$ and $f_{-t}\ dh_{\gamma}^*(\nu)=
\nu$.

For any $t\in {\RR}$, denote by $Z_t$ the set of relative fixed points
of $f_t$.
We also introduce the corresponding set in the cospherical bundle
$SN^*{\cal F}=\{\nu\in N^*{\cal F}: \sigma_P(\nu)=1\}$:
$S{Z}_{t}={Z}_{t}\cap SN^*{\cal F}.$
This set might be not closed, but, for any
$k\in C^{\infty}_c(G,|T{\cG}|^{1/2}))$,
the corresponding part
$$S{Z}_{t,k}= \{\nu\in SN^*{\cal F} : (\exists \gamma \in
{\rm supp}\,k, r(\gamma)=\pi(\nu)) f_{-t}\ dh_{\gamma}^*(\nu)=\nu\}$$
is closed. By the tranversal ellipticity of $\sigma_A$, the flow $f_t$ is
transverse to ${\cF}_N$, therefore, {\bf the relative period set}
${\cT}_k=\{t\in {\RR}: S{Z}_{t,k} \not= \varnothing   \}$
is a discrete subset of $\RR$.

The following theorem was proved in \cite{geo}, but we will give
an independent proof.

\begin{thm}
\label{relation}
Given an operator $A$ under the conditions (A1) and (A2) and $k\in
C^{\infty}_c(G,|T\cG|^{1/2})$, the distribution $\theta_k$
is smooth outside of the relative period set ${\cal T}_k$ of the transverse
bicharacterictic flow $f_t$.
\end{thm}

Let $G_{{\cal F}_N}$ denote the holonomy groupoid of the foliation
$\cF_N$. $G_{{\cal F}_N}$ consists of all pairs
$(\gamma,\nu)\in G_{\cal F}\times \tilde{N}^{*}{\cal F}$ such that
$r(\gamma)=\pi(\nu)$
with the source map $s_N:G_{{\cal F}_N}\rightarrow
\tilde{N}^{*}{\cF}, s_N(\gamma,\nu)=dh_{\gamma}^{*}(\nu)$, and the target map
$r_N:G_{{\cal F}_N}\rightarrow \tilde{N}^{*}{\cal F}, r_N(\gamma,\nu)=\nu$.
There is a projection $\pi_G:G_{{\cal F}_N}\rightarrow G_{\cal F}$
given by the formula $\pi_G(\gamma,\nu)=\gamma$.
Put also $G_{SN^*{\cF}}=G_{{\cal F}_N}\cap r^{-1}_N(SN^*{\cal F})$.

For any $(\gamma,\nu)\in G_{{\cF}_N}$, denote by
$dH_{(\gamma,\nu)}$ the associated linear holonomy map:
$$
dH_{(\gamma,\nu)}:
H_{dh_{\gamma}^*(\nu)}{{\cal F}_N}\rightarrow H_{\nu}{{\cal F}_N}.
$$
It is easy to see that $dH_{(\gamma,\nu)}$ preserves
the symplectic structure of $H{{\cal F}_N}$.

Denote by $Q:TN^*\cF\to H\cF_N$ the projection map.
The differential of the map $(r_N,s_N): G_{{\cal F}_N}\to \tilde{N}^*{\cal F}
\times\tilde{N}^*{\cal F}$ defines an isomorphism of the tangent
space $T_{(\gamma,\nu)}G_{{\cal F}_N}$ with the set of all
$(v_1,v_2)\in T_{\nu}N^*{\cF}\oplus T_{dh^*_{\gamma}(\nu)}N^*{\cF}$ such that
the normal components of $v_1$ and $v_2$ are connected by the holonomy map:
$Q(v_1)=dH_{(\gamma,\nu)}(Q(v_2))$.

The holonomy groupoid $G_{{\cal F}_N}$ has the natural foliation
${\cG}_{{\cal F}_N}$ such that the tangent bundle
$T_{(\gamma,\nu)}{\cG}_{{\cal F}_N}$ corresponds to
$T_{\nu}{\cF}_N\oplus T_{dh^*_{\gamma}(\nu)}{\cF_N}$ under the isomorphism
described above. The normal space $H_{(\gamma,\nu)}{{\cG}_{{\cal F}_N}}$
to $\cG_{{\cal F}_N}$  is isomorphic to the set of all
$(v_1,v_2)\in H_{\nu}{\cF}_N\oplus H_{dh^*_{\gamma}(\nu)}{\cF}_N$ such that
$v_1=dH_{(\gamma,\nu)}(v_2)$, and therefore the maps $dr_N$ ($ds_N$)
define isomorphisms of $H_{(\gamma,\nu)}{{\cG}_{{\cal F}_N}}$ with
$H_{\nu}{{\cal F}_N}$ ($H_{dh^*_{\gamma}(\nu)}{\cF}_N$) respectively.

\begin{lem}\label{l:saturated}
The set $Z_t$ is a saturated subset of $N^*{\cal F}$, that is, it is a
union of leaves
of the foliation $\cF_N$.
\end{lem}
\begin{proof} By the holonomy invariance of $p$, the Hamiltonian vector
field $\Xi_p$
with the Hamiltonian $p$ satisfies the identity
$$
dH_{(\gamma,\nu)}(Q(\Xi_p(dh^*_{\gamma}(\nu))))=Q(\Xi_p(\nu)),\quad
\nu\in \tilde{N}^*{\cF},\quad \gamma\in G,\quad r(\gamma)=\pi(\nu),
$$
therefore, there exists a vector field $\hat{\Xi}_p$ on $G_{{\cal F}_N}$
such that
\begin{equation}\label{e:lift}
ds_N(\hat{\Xi}_p(\gamma,\nu))=\Xi_p(dh^*_{\gamma}(\nu)),\quad
dr_N(\hat{\Xi}_p(\gamma,\nu))=\Xi_p(\nu), \quad (\gamma,\nu)\in G_{{\cal F}_N}.
\end{equation}

Let $\hat{F}_t$ be a flow on $G_{{\cal F}_N}$ generated by $\hat{\Xi}_p$.
By \eqref{e:lift}, we have
$$
f_t\circ r_N=r_N\circ \hat{F}_t,\quad f_t\circ s_N=s_N\circ \hat{F}_t,
$$
or, if we write $\hat{F}_t:G_{{\cal F}_N}\to G_{{\cal F}_N}$ as
$\hat{F}_t(\gamma,\nu)=(F_t(\gamma,\nu),f_t(\nu))$,
\begin{equation}\label{e:shift}
f_t(dh_{\gamma}^*(\nu))=dh_{F_t(\gamma,\nu)}^*(f_t(\nu)).
\end{equation}

Take any $\nu\in Z_t$ with the corresponding $\gamma\in G$ such
$r(\gamma)=\pi(\nu), f_{-t}\ dh_{\gamma}^*(\nu)=\nu$. Let
$(\gamma_1,\nu)\in G_{{\cal F}_N}$. Then we have
\begin{alignat*}{2}
f_{t}(dh_{\gamma_1}^{*}(\nu))&=dh^*_{F_t(\gamma_1,\nu)}(f_t(\nu))\qquad
\text{by \eqref{e:shift}}\\
&=dh^*_{F_t(\gamma_1,\nu)}(dh_{\gamma}^*(\nu))\\
&=(dh^*_{F_t(\gamma_1,\nu)}\circ dh_{\gamma}^*\circ dh_{\gamma_1^{-1}}^{*})
(dh_{\gamma_1}^{*}(\nu))\\ &=dh_{{\gamma'}}^{*}(dh_{\gamma_1}^{*}(\nu)),
\end{alignat*}
where $\gamma'=\gamma^{-1}_1\circ \gamma \circ F_t(\gamma_1,\nu)$ that implies
$dh_{\gamma_1}^{*}(\nu)\in Z_t$.
 \end{proof}

The relative fixed point sets $Z_t$ can be naturally lifted to the holonomy
groupoid $G_{{\cal F}_N}$:
\begin{equation*}
{\cal Z}_{t}= \{(\gamma,\nu)\in G_{{\cal F}_N}:
f_{-t}\ dh_{\gamma}^*(\nu)=\nu\}, \quad
S{\cal Z}_{t}={\cal Z}_{t}\cap G_{SN^*{\cal F}},
\end{equation*}
By Lemma~\ref{l:saturated}, ${Z}_{t}=r_N({\cal Z}_{t})=s_N({\cal Z}_{t})$.

Let us assume that ${\cal Z}_{t}$ is a smooth submanifold
of $G_{{\cal F}_N}$. By Lemma~\ref{l:saturated}, the tangent space to $\cZ_t$
at a point $(\gamma,\nu)\in \cZ_t$ contains a subspace
$F_{(\gamma,\nu)}{\cZ_t}$,
which is the graph of the linear map
$df_t(\nu):T_{\nu}{\cF}_N\to T_{dh^*_\gamma(\nu)}{\cF}_N=T_{f_t(\nu)}{\cF}_N$:
$$
F_{(\gamma,\nu)}{\cZ_t}=\{(v_1,v_2)\in T_{\nu}{\cF}_N\times
T_{dh^*_\gamma(\nu)}{\cF}_N: v_2=df_t(\nu)(v_1) \}.
$$
Let
$$
H_{(\gamma,\nu)}{\cZ_t}=T_{(\gamma,\nu)}\cZ_t/F_{(\gamma,\nu)}{\cZ_t}, \quad
H_{(\gamma,\nu)}{S\cZ_t}=T_{(\gamma,\nu)}S\cZ_t/F_{(\gamma,\nu)}{\cZ_t}.
$$

\begin{defn}\label{clean}
Let $t\in\RR$ be a relative period of the flow $f_t$. We say that
the flow $f_t$ is {\bf clean} on ${\cal Z}_{t}$, if:
\begin{itemize}
\item[(1)] ${\cal Z}_{t}$ is a smooth submanifold
of $G_{{\cal F}_N}$;

\item[(2)] the normal space $H_{(\gamma,\nu)}{\cZ_t}$
at any point $(\gamma,\nu)\in \cZ_t$ coincides with the set of
all $(v_1,v_2)\in H_{(\gamma,\nu)}{{\cG}_{{\cal F}_N}}$ such that
$v_2=df_t(\nu)(v_1).$
\end{itemize}
\end{defn}

Let $|T{\cal F}_N|^{1/2}$ be
the vector bundle of leafwise half-densities on $N^*{\cal F}$, and
$s_N^*(|T{\cal F}_N|^{1/2})$ and $r_N^*(|T{\cal F}_N|^{1/2})$
are the lifts of this vector bundle to vector bundles on $G_{{\cal F}_N}$
via the mappings $s_N$ and
$r_N$ respectively.
Let $|T{\cG}_{{\cal F}_N}|^{1/2}$ be the vector bundle of
leafwise half-densities on $G_{{\cal F}_N}$:
$$
|T{\cG}_{{\cal F}_N}|^{1/2}=r_N^*(|T{\cal F}_N|^{1/2})\otimes
s_N^*(|T{\cal F}_N|^{1/2}).
$$
The projection $\pi_G:G_{\cF_N}\to G$ defines a local diffeomorphism
$\pi_G:\cG_{\cF_N}\to \cG$, that induces a map
$$
\pi^*_G:C^{\infty}_c(G,|T\cG|^{1/2})\to
C^{\infty}(G_{\cF_N},|T\cG_{\cF_N}|^{1/2}).
$$

Define a restriction map
$$
R_\cZ:C^{\infty}_c(G_{\cF_N},|T\cG_{\cF_N}|^{1/2})\to
C^{\infty}_c(\cZ_t,|T\cF_N|^1)
$$
as follows. If $\rho=fr_N^{\ast}\rho_1\otimes
s_N^{\ast}\rho_2, \quad f\in C^{\infty}_c(G_{\cF_N}), \quad
\rho_1,\rho_2\in C^{\infty}_c(M,|T\cF_N|^{1/2})$,
then
$$
R_\cZ\rho(\gamma,\nu)=f(\gamma,\nu)\rho_1(\gamma,\nu)
df^*_{t}(\nu)[\rho_2(dh^*_{\gamma}(\nu))],\quad (\gamma,\nu)\in {\cZ}_t,
$$
where the map
$df_{t}^*(\nu): |T_{f_t(\nu)}{\cal F}_N|^{1/2}\rightarrow
|T_{\nu}{\cal F}_N|^{1/2}$
is induced by the linear map
$df_t(\nu):T_{\nu}{\cF}_N\to T_{f_t(\nu)}{\cF}_N$.

If the flow $f_t$ is clean, there is defined a natural density $d\mu_{\cZ}$
on $H_{(\gamma,\nu)}\cZ_t$, being the fixed point set  of
the symplectic linear map $dH_{(\gamma,\nu)}\circ df_t(\nu)$
of the symplectic space $H_{\nu}{\cal F}_N$
(see, for instance, \cite[Lemma 4.3]{DG}).
Dividing $d\mu_{\cZ}$ by $d\sigma_P$, we get a density
$d\mu_{S\cZ}$ on $H_{(\gamma,\nu)}S\cZ_t$.

Using the natural isomorphism
\begin{equation*}
|TS\cZ_t|\cong |F{\cZ}| \otimes |HS\cZ_t|.
\end{equation*}
one can combine the densities $R_\cZ\pi^*_Gk\in
C^{\infty}_c(S\cZ_t, |FS\cZ_t|)$ and $d\mu_{S\cZ}\in
C^{\infty}_c(S\cZ_t,|HS\cZ_t|)$ to
get a smooth density $R_\cZ\pi^*_Gk\,d\mu_{S{\cZ}}$ on $S{\cZ}_t$.

Let $\sigma_{\rm sub}(A)$ denote the subprincipal symbol of $A$. Define
$\sigma_{\rm sub}(P)=\frac{1}{2}a^{-\frac{1}{2}}\sigma_{\rm sub}(A)$ in some
conic neighborhood of $N^*{\cF}$. The restriction of $\sigma_{\rm sub}(P)$
on $N^*{\cF}$
is equal to the restriction on $N^*{\cF}$ of the subprincipal symbol of
any operator $P_1\in \Psi^1(M,|TM|^{1/2})$
such that the complete symbols of $P^2_1$ and $A$ are equal
$\mod S^{-\infty}$ in some neighborhood of $N^*{\cal F}$.

\begin{thm}\label{t:main}
Let $t\in\RR$ be a relative period of the flow $f_t$. Assume that the
relative fixed point set
${\cal Z}_t$ is clean.
Then, for any $k\in C^{\infty}_c(G,|T\cG|^{1/2})$ and for any $\tau $
in some neighborhood of $t$, we have
\begin{equation}
\label{Ip}
\theta_k(\tau)=\sum_{{\cal Z}_j} \int_{-\infty}^{+\infty} \alpha_j(s,k)
e^{is(\tau-t)}ds,
\end{equation}
where:
\begin{enumerate}
\item ${\cal Z}_j$ are all connected components of the set
$S{\cal Z}_t$ in $G_{SN^*{\cal F}}$ of dimensions $d_j=\dim {\cal Z}_j$;
\item $\alpha_j$ has an asymptotic expansion
\begin{equation}\label{expansion}
\alpha_j(s,k)\sim \left(\frac{s}{2\pi i}\right)^{(d_j-p-1)/2} i^{-\sigma_j}
\sum_{r=0}^{+\infty}\alpha_{j,r}(k)s^{-r}, \quad s\rightarrow +\infty
\end{equation}
with $\alpha_{j,0}$ given by the formula
\begin{equation}\label{e:leading}
\alpha_{j,0}(k)  = \int_{{\cal Z}_j}
e^{i\int_0^t \sigma_{\rm sub}(P)(f_{-\tau}dh_{\gamma}^*(\nu))d\tau}
R_\cZ\pi^*_Gk(\gamma,\nu) d\mu_{{S\cZ}_j}(\gamma,\nu),
\end{equation}
where $\sigma_j$ denotes the Maslov index associated with the connected
component ${\cal Z}_j$ (see below for the definition).
\end{enumerate}
\end{thm}

\section{Reduction to the case when $A$ is elliptic}
In this section, we will assume that $A$ is an operator under the assumptions
(A1) and (A2).
We will use the classes $\Psi^{m,-\infty}(M,{\cal F},|TM|^{1/2})$ of
transversal
pseudodifferential operators (see \cite{noncom} for the definition) and
the Sobolev spaces $H^s(M)$ of half-densities on $M$.
Put also $\Psi^{*,-\infty}(M,{\cal F},|TM|^{1/2})=
\bigcup_m\Psi^{m,-\infty}(M,{\cal F},|TM|^{1/2})$.

By \cite{noncom}, the operator $P=A^{1/2}$ satisfies
the following conditions:
\bigskip\par
(H1) $P$ has the form
$$P=P_1+R_1,$$
where:
\medskip\par
(a) $P_1\in \Psi^1(M,|TM|^{1/2})$ is a transversally elliptic operator with
the positive, holonomy invariant transversal principal symbol such that
the complete symbols of $P^2_1$ and $A$ are equal $\mod S^{-\infty}$
in some neighborhood of $N^*{\cal F}$;
\medskip\par
(b) $R_1$ is a bounded operator from $L^2(M)$ to $H^{-1}(M)$ and
for any $K\in \Psi^{*,-\infty}(M,{\cal F},|TM|^{1/2})$ the operator $KR_1$ is
a smoothing operator in $L^2(M)$, that is, it defines a bounded operator from
$L^2(M)$ to $C^{\infty}(M,|TM|^{1/2})$.
\medskip\par
(H2) $P$ is essentially self-adjoint in $L^2(M)$ (with the initial domain
$C^{\infty}(M,|TM|^{1/2})$).

\begin{lemma}[\cite{geo}]
Any operator $P$, satisfying the conditions (H1) and (H2),
can be represented in the form
\begin{equation}
\label{elliptic}
P=P_2+R_2,
\end{equation}
where:
\medskip\par
(a) $P_2\in \Psi^1(M,|TM|^{1/2})$ is an essentially self-adjoint,
elliptic operator with the
positive principal symbol and the holonomy invariant transversal
principal symbol such that the complete symbols of $P_1$ and $P_2$
are equal $\mod S^{-\infty}$ in some neighborhood of $N^*{\cal F}$;
\medskip\par
(b) $R_2$ is a bounded operator from $L^2(M)$ to $H^{-1}(M)$ and,
for any $K\in \Psi^{*,-\infty}(M,{\cal F},|TM|^{1/2})$, the operator $KR_2$ is
a smoothing operator in $L^2(M)$.
\end{lemma}

\begin{proof}
Take a foliated coordinate chart $\Omega$ on $M$ with coordinates
$(x,y)\in I^p\times I^q$ ($I$ is the open interval $(0,1)$)
such that the restriction of $\cF$ on $U$ is given by the sets
$y={\rm const}$.
Let $p_1\in S^1(I^n\times {\RR}^n)$ be the
complete symbol of the operator $P_1$ in this chart.
Assume that $p_1(x,y,\xi,\eta)$ is invertible for any
$(x,y,\xi,\eta)\in U, |\xi|^2+|\eta|^2>R^2$, where $R>0$, $U$ is a
conic neighborhood of the set $\eta=0$.
Take any function $\phi\in C^{\infty}(I^n\times {\RR}^n),
\phi=\phi(x,y,\xi,\eta), x\in I^p, y\in I^q, \xi \in {\RR}^p,
\eta \in {\RR}^q$, homogeneous of degree $0$ in $(\xi,\eta)$
for $|\xi|^2+|\eta|^2>1$, which is supported in some conic neighborhood
of $\eta=0$ and is equal to $1$ in $U$, and put
$$
p_2(x,y,\xi,\eta)=\phi p_1(x,y,\xi,\eta) + (1-\phi)
(1+|\xi|^2+|\eta|^2)^{1/2}.
$$
Take $P_2$ to be the operator $p_2(x,y,D_x,D_y)$
with the complete symbol $p_2$ (or, more precisely,
$p_2(x,y,D_x,D_y)+p_2(x,y,D_x,D_y)^*$ to provide self-adjointness)
and put $R_2=P-P_2$. The operator $P_1-P_2$ has order $-\infty$ in some
conic neighbourhood of $N^*{\cF}$, therefore, for any
$K\in \Psi^{*,-\infty}(M,{\cal F},|TM|^{1/2})$ the operator $K(P_1-P_2)$
is a smoothing operator
\cite{noncom}, that completes immediately the proof.
\end{proof}

Denote by $W(t)=e^{itP_2}$ the wave group generated by the elliptic
operator $P_2$.
It is well-known that $W(t)$ is a Fourier integral operator (see below for
more details).
Put also $R(t)=e^{itP}-W(t)$.

\begin{proposition}
\label{mapping}
For any $K\in \Psi^{*,-\infty}(M,{\cal F},|TM|^{1/2})$, the family
$KR(t), t\in \RR,$ is
a smooth family of bounded operators from $L^2(M)$ to
$C^{\infty}(M,|TM|^{1/2})$.
\end{proposition}

\begin{proof} Since $P^2=A\in \Psi^2(M,|TM|^{1/2})$,
by interpolation and duality, $P$ defines a bounded operator
from $H^1(M)$ to $L^2(M)$ and from $L^2(M)$ to $H^{-1}(M)$ and, for any
natural $N$, $P^N$ defines a bounded operator from $H^N(M)$ to $L^{2}(M)$
and from $L^2(M)$ to $H^{-N}(M)$.
Since $R_2=P-P_2$ and $P_2\in \Psi^1(M,|TM|^{1/2})$, $R_2$ also defines
a bounded operator from $H^1(M)$ to $L^2(M)$ and from $L^2(M)$ to $H^{-1}(M)$.

By assumption,  for any $K\in \Psi^{*,-\infty}(M,{\cal F},|TM|^{1/2})$,
the operator $KR_2$ is
a smoothing operator, therefore, the operator $KR_2P^N$ is defined as an
operator from $H^N(M)$ to $C^{\infty}(M,|TM|^{1/2})$.

\begin{lem}\label{l:power}
For any $K\in \Psi^{*,-\infty}(M,{\cal F},|TM|^{1/2})$ and for any $N\in\NN$,
the operator $KR_2P^N$ extends to a bounded operator from $L^2(M)$
to $C^{\infty}(M,|TM|^{1/2})$.
\end{lem}
\begin{proof}
We will prove the lemma by induction on $N$. For $N=0$, the statement is
true by assumption. Let us assume that it is true for some $N$, that is,
for any $K\in \Psi^{*,-\infty}(M,{\cal F},|TM|^{1/2})$,
the operator $KR_2P^N$ extends to a bounded operator from $L^2(M)$
to $C^{\infty}(M,|TM|^{1/2})$.

We have the equality $P^2-P^2_2=R_2P+P_2R_2$ as operators from $H^1(M)$
to $H^{-1}(M)$,
therefore, for any $K\in \Psi^{*,-\infty}(M,{\cal F},|TM|^{1/2})$,
\begin{equation*}
KR_2P^{N+1}=K(R_2P)P^N=K(P^2-P^2_2)P^N-KP_2R_2P^N.
\end{equation*}
The operator $P^2-P^2_2\in \Psi^{2}(M,|TM|^{1/2})$ has order $-\infty$ in some
conic neighbourhood of $N^*{\cF}$, therefore, for any
$K\in \Psi^{*,-\infty}(M,{\cal F},|TM|^{1/2})$ the operator $K(P^2-P^2_2)$
extends to a
bounded operator from $H^s(M)$ to $C^{\infty}(M,|TM|^{1/2})$ for any $s$ and
$K(P^2-P^2_2)P^N$ extends to a
bounded operator from $L^2(M)$ to $C^{\infty}(M,|TM|^{1/2})$.

Since $P_2\in \Psi^{1}(M,|TM|^{1/2})$ and
$K\in \Psi^{*,-\infty}(M,{\cal F},|TM|^{1/2})$,
by the composition theorem \cite{noncom},
$KP_2\in  \Psi^{*,-\infty}(M,{\cal F},|TM|^{1/2})$
and, by induction hypothesis, $KP_2R_2P^N$ extends to a
bounded operator from $L^2(M)$ to $C^{\infty}(M,|TM|^{1/2})$.
\end{proof}

By the Duhamel formula, we have
\begin{equation*}
R(t)u=i \int_0^t e^{i\tau P_2}\,R_2\,e^{i(t-\tau)P}u\,d\tau,
\quad u\in H^1(M)\subset D(P),
\end{equation*}
therefore, for any $K\in \Psi^{*,-\infty}(M,{\cal F},|TM|^{1/2})$,
\begin{equation*}
KR(t)=i \int_0^t Ke^{i\tau P_2}\,R_2\,e^{i(t-\tau)P}\,d\tau
         =i \int_0^t e^{i\tau P_2}e^{-i\tau P_2}Ke^{i\tau P_2}\,R_2\,
e^{i(t-\tau)P}\,d\tau.
\end{equation*}
Any operator $K\in \Psi^{*,-\infty}(M,{\cal F},|TM|^{1/2})$ is
a Fourier integral operator
(see below for more details) and, using the composition theorem for
Fourier integral operators, one can check that
$e^{-i\tau P_2}Ke^{i\tau P_2}\in \Psi^{*,-\infty}(M,{\cal F},|TM|^{1/2})$.
Therefore, the
operator $e^{-i\tau P_2}Ke^{i\tau P_2}R_2$ extends to a
bounded operator from $L^2(M)$ to $C^{\infty}(M,|TM|^{1/2})$. Since
$e^{i\tau P_2}$ maps $C^{\infty}(M,|TM|^{1/2})$ to
$C^{\infty}(M,|TM|^{1/2})$ and,
by the spectral theorem, $e^{i(t-\tau)P}$ is a bounded operator in $L^2(M)$,
for any $K\in \Psi^{*,-\infty}(M,{\cal F},|TM|^{1/2})$, the operator $KR(t)$
extends to a bounded operator from $L^2(M)$ to $C^{\infty}(M,|TM|^{1/2})$.
Moreover, one can be easily seen from above arguments that the function
$KR(t)$ is continuous as a function on $\RR$ with values in the space
$\cL(L^2(M),C^{\infty}(M,|TM|^{1/2}))$
of bounded operators from $L^2(M)$ to $C^{\infty}(M,|TM|^{1/2})$.

For any $u\in H^1(M)$, the function $\RR\to H:t\mapsto KR(t)u$ is
differentiable,
and
$$
\frac{d}{dt}KR(t)u=iK(Pe^{itP}u-P_2e^{itP_2}u)=i(KP_2R(t)+KR_2e^{itP})u.
$$
The operator $KP_2R(t)+KR_2e^{itP}$ extends to a bounded operator
from $L^2(M)$ to $C^{\infty}(M,|TM|^{1/2})$, and, moreover, the function
$t\mapsto KP_2R(t)+KR_2e^{itP}$ is a continuous function on $\RR$ with
values in $\cL(L^2(M),C^{\infty}(M,|TM|^{1/2}))$.
Using this, one can be easily seen that the function
$t\mapsto KR(t)$ is differentiable
as a function on $\RR$ with values in $\cL(L^2(M),C^{\infty}(M,|TM|^{1/2}))$
and
$$
\frac{d}{dt}KR(t)=i(KP_2R(t)+KR_2e^{itP}).
$$

Let us proceed by induction.  Assume that, for any
$K\in \Psi^{*,-\infty}(M,{\cal F},|TM|^{1/2})$ and for some natural $n$,
the function $KR(t)$ is $n$-times differentiable
as a function on $\RR$ with values in $\cL(L^2(M),C^{\infty}(M,|TM|^{1/2}))$
and the derivative $KR^{(n)}(t), t\in \RR$
satisfies the equation
\begin{equation}\label{e:induction0}
KR^{(n)}(t)=iKP_2R^{(n-1)}(t)+i^{n}KR_2P^{n-1}e^{itP}.
\end{equation}

To prove that the function $t\mapsto KR^{(n)}(t)$ is differentiable
as a function on $\RR$ with values in $\cL(L^2(M),C^{\infty}(M,|TM|^{1/2}))$,
as above, it suffices to prove that the derivative $(d/dt)KR^{(n)}(t)u$
exists for any $u$ from a dense subspace of $L^2(M)$, it extends to a bounded
operator from $L^2(M)$ to $C^{\infty}(M,|TM|^{1/2})$, and its
extension is continuous as a function on $\RR$ with values
in $\cL(L^2(M),C^{\infty}(M,|TM|^{1/2}))$.

From \eqref{e:induction0}, one can easy to see that the derivative
$(d/dt)KR^{(n)}(t)u$
exists for any $u\in H^1(M)$ and satisfies the equation
\begin{equation}\label{e:induction1}
\frac{d}{dt}KR^{(n)}(t)u=iKP_2R^{(n)}(t)u+i^{n+1}KR_2P^ne^{itP}u.
\end{equation}
The first term in the right-hand side of \eqref{e:induction1},
$iKP_2R^{(n)}(t)$,
is a  bounded operator from $L^2(M)$ to $C^{\infty}(M,|TM|^{1/2})$ by the
induction hypothesis. By Lemma~\ref{l:power}, the operator $KR_2P^n$
extends to a bounded operator from $L^2(M)$ to $C^{\infty}(M,|TM|^{1/2})$,
and, by the spectral theorem, $e^{itP}$ is a bounded operator in $L^2(M)$,
therefore,
the second term in the right-hand side of \eqref{e:induction1},
the operator $i^{n+1}KR_2P^ne^{itP}$, extends to a bounded operator
from $L^2(M)$ to $C^{\infty}(M,|TM|^{1/2})$. It is also clear the
right-hand side of \eqref{e:induction1} is continuous as a function
on $\RR$ with values in $\cL(L^2(M),C^{\infty}(M,|TM|^{1/2}))$.
This completes the proof of the existence of the derivative $KR^{(n+1)}(t)=
(d/dt)KR^{(n)}(t)$ and the induction arguments.
\end{proof}

\begin{proof}[Proof of Proposition~\ref{p:existence}]
Let $W(t)$ and $R(t)$ be as in Proposition~\ref{mapping} and
$k\in C_c^{\infty}(G,|T\cG|^{1/2})$. Define $\theta_k(t)$
by the formula
$$
\theta_k(t)={\rm tr}\,R(k)W(t)+{\rm tr}\,R(k)R(t).
$$
Since $P_2$ is an elliptic operator, the operator $\int f(t)e^{itP_2}dt$
is a smoothing operator in $\cD'(M)$ and the trace of the
operator $R(k)W(t)$ is well-defined
as a distribution on $\RR$ \cite{DG}.
Since any bounded operator $T$ in $L^2(M)$, which extends to a bounded
operator from $L^2(M)$ to $H^s(M)$ with $s>n=\dim M$, is a trace class
operator, the trace
${\rm tr}\, R(k)R(t)$ is a well-defined smooth function on ${\RR}$ by
Proposition~\ref{mapping}.
\end{proof}

\begin{cor}
For any $k\in C_c^{\infty}(G,|T\cG|^{1/2})$, the function
${\rm tr}\,R(k)R(t)=\theta_k(t)-{\rm tr}\,R(k)W(t)$ is a smooth function
on $\RR$.
\end{cor}

It should be noted that, without any additional assumption about the
operator $P$ in question, the
corresponding distribution on $G_{{\cal F}_N}\times {\RR}$,
$k \mapsto {\rm tr}\, R(k)R(t)$, might be very singular, but
the singularities of the distribution $k \mapsto {\rm tr}\, R(k)W(t)$
can be described rather explicitly under the clear
intersection assumption.

\section{The case of an elliptic operator}
Let $P_2\in \Psi^1(M,|TM|^{1/2})$ be an essentially self-adjoint,
elliptic operator with the positive principal symbol and the holonomy
invariant transversal
principal symbol and $W(t)=e^{itP_2}$.
The singularities of the distribution $t\mapsto {\rm tr}\, R(k)W(t)$
can be studied in a standard manner, using microlocal analysis.

Fix $k\in C^{\infty}_c(G,|T\cG|^{1/2})$. We will consider the operator
family  $R(k)W(t)$
as a single operator $R(k)W$ from $L^2(M)$ to $L^2({\RR}\times M)$.
We will prove that this operator is a Fourier integral operator. At first,
let us recall well-known facts about the structure of
the operators $R(k)$ and $W$.

As above, let ${\tilde p}$ be a smooth function on $\tilde{T}^*M$
homogeneous of degree one such that $\tilde{p}(\xi)\not=0$ for
$\xi\in \tilde{T}^*M$, which is equal to $p=a^{1/2}$ in some conic neighborhood
of $N^*{\cal F}$ and  $\tilde{f}_t$ the Hamiltonian flow of $\tilde{p}$.
Without loss of generality, we may assume that $\tilde{p}$ is the principal
symbol of the operator $P_2$. Let $\Lambda_{\tilde p}$ be the Lagrangian
submanifold
in $\tilde{T}^*{\RR}\times \tilde{T}^*M \times \tilde{T}^*M$:
\begin{equation*}
\Lambda_{\tilde p}= \{((t,\tau),(x,\xi),(y,\eta))\in
\tilde{T}^*{\RR}\times\tilde{T}^*M\times \tilde{T}^*M:
\tau = \tilde{p}(x,\xi), (x,\xi)=\tilde{f}_{-t}(y,\eta)\}.
\end{equation*}
Then $W$ is a Fourier integral operator associated with the canonical
relation $\Lambda'_{\tilde p}$, $W\in I^{-1/4}({\RR}\times M\times
M,\Lambda'_{\tilde p})$.

The operator $R(k)$ belongs to $\Psi^{0,-\infty}(M,{\cal F},|TM|^{1/2})$  and,
therefore, is a Fourier integral operator associated with an immersed
canonical relation, which is the image of $G_{{\cal F}_N}$ under
the mapping
$$
(r_N,s_N):G_{{\cal F}_N}\rightarrow T^*M\times T^*M,\quad (\gamma,\nu)\mapsto
(\nu, dh_{\gamma}^{*}(\nu)),
$$
given by the source and the target mappings of
the groupoid $G_{{\cal F}_N}$ (see \cite{noncom}).
More precisely, $R(k)\in I^{-p/2}(M\times M,G'_{{\cal F}_N})$.

By the tranversal ellipticity of $\tilde{p}$, the intersection of
$\Lambda'_{\tilde p}$ with $G'_{{\cal F}_N}$ is transverse, and,
by the composition theorem of Fourier integral
operators \cite{H4},
the operator $R(k)W$ is a Fourier integral operator associated with an
immersed canonical relation from $T^*M$ to $T^*({\RR}\times M)$
given by the map
$$
\Pi:{\RR}\times G_{{\cal F}_N}\rightarrow T^*{\RR}\times
T^*M\times T^*M, \quad
(t,\gamma,\nu)\mapsto (t,p(\nu),\nu,f_{-t}dh^*_{\gamma}(\nu)).$$
More precisely, $R(k)W\in I^{-p/2-1/4}({\RR}\times M\times M;
{\RR}\times G_{{\cal F}_N},\Pi).$

Recall that the trace functional can be treated from the point
of microlocal analysis as follows \cite{DG}.
Let $\Delta:{\RR}\times M\rightarrow {\RR}\times M \times M$ be
the diagonal map, $\Delta(t,x)=(t,x,x), (t,x)\in {\RR}\times M$,
and $\pi : {\RR}\times M\rightarrow M$ the projection
map. Then
\begin{equation}
\label{trace}
\hbox{tr}\ R(k)W= \pi_*\Delta^*W_k,
\end{equation}
where
$
W_k\in C^{\infty}({\RR}\times M \times M,
|T({\RR}\times M \times M)|^{1/2})
$
is the Schwartz kernel of the operator $R(k)W$,
$\Delta^*:C^{\infty}({\RR}\times M \times M,
|T({\RR}\times M \times M)|^{1/2})\rightarrow
C^{\infty}({\RR},|T{\RR}|^{1/2})\otimes C^{\infty}(M,|TM|)
$
is defined by the formula
$$
\Delta^*(s_1\otimes s_2\otimes s_3)(t,x)=s_1(t)\otimes (s_2(x)\otimes
s_3(x)),\quad  t\in {\RR}, \quad x\in M,
$$
where $s_1\in C^{\infty}({\RR}, |T{\RR}|^{1/2})$,
$s_2\in C^{\infty}(M,|TM|^{1/2})$,
$s_3\in C^{\infty}(M,|TM|^{1/2})$,
and
$$
\pi_*:C^{\infty}({\RR},|T{\RR}|^{1/2})\otimes
C^{\infty}(M,|TM|)\rightarrow
C^{\infty}({\RR},|T{\RR}|^{1/2})
$$
is given by integration along fibers of the projection $\pi$.

It is known that $\pi_*\Delta^*\in I^0({\RR}\times M\times M\times
{\RR},\Gamma)$, where $\Gamma$ is the conormal bundle to the
diagonal in ${\RR}\times M\times M\times {\RR}$:
$$
\Gamma=\{(t,\tau_1,\nu_1,\nu_2,t,\tau_2)\in
T^*{\RR}\times T^*M\times T^*M\times T^*{\RR}:
\nu_1=-\nu_2, \tau_1=-\tau_2\}.
$$
There is a commutative diagram
\begin{equation}\label{diag}
\CD
\Gamma @<p_1<< {\cZ}\\
@V\varphi VV @Vp_2 VV\\
T^*({\RR}\times M\times M) @<\Pi << {\RR}\times G_{{\cal F}_N}
\endCD
\end{equation}
where
$$
p_1(t,\gamma,\nu) = (\Pi(t,\gamma,\nu),t,-p(\nu))=
(t, p(\nu), \nu, -\nu, t, -p(\nu)),\quad  (t,\gamma,\nu)\in {\cZ},
$$
$p_2$ is a natural inclusion and
$$
\varphi(t,\tau,\nu,-\nu,t,-\tau)=(t,\tau,\nu,-\nu), \quad
(t, \tau, \nu, -\nu, t, -\tau) \in \Gamma.
$$
It is easy to see that (\ref{diag}) is a fiber product diagram, that is,
$$
{\cZ} \cong \{(x,y)\in ({\RR}\times G_{{\cal F}_N})\times
 \Gamma : \Pi(x)=\varphi(y)\}.
$$
Using this fact and the functoriality properties of the
wave-front sets (see, for instance, \cite{geom:asymp}),
one get immediately the description of the
singularities of the distribution $\theta_k$,
given by Theorem~\ref{relation}.
\medskip\par
To finish the proof of Theorem~\ref{t:main},
we will state under the assumption on the flow $f_t$
to be clean in the sense of
Definition~\ref{clean} that $\pi_*\Delta^*W_k$ is
a Lagrangian distribution and compute its symbol.
We begin with computation of the symbol of the operator $R(k)W$.
Recall first the description of the principal symbol of the
operator $R(k)$.

According to the short exact sequence
$$
 0\rightarrow T\cG_{\cF_N}\rightarrow TG_{\cF_N}\rightarrow H\cG_{\cF_N}
\rightarrow 0,
$$
the half-density vector bundle on $G_{{\cal F}_N}$ can be decomposed as
$$
|TG_{{\cal F}_N}|^{1/2}\cong |T{\cG}_{{\cal F}_N}|^{1/2}
\otimes |H{{\cG}_{{\cal F}_N}}|^{1/2},
$$
where $|H{{\cG}_{{\cal F}_N}}|^{1/2}$ is
the transverse half-density bundle on $G_{{\cal F}_N}$:
$|H_{(\gamma,\nu)}{{\cG}_{{\cal F}_N}}|^{1/2}
\cong |H_{\nu}{{\cal F}_N}|^{1/2}
 \cong |H_{dh^*_{\gamma}(\nu)}{{\cal F}_N}|^{1/2}.$

Let $|dy\wedge d\eta|^{1/2}\in C^{\infty}(N^*{\cal F},
|H{{\cal F}_N}|^{1/2})$ be given by the Liouville form of
the canonical transverse symplectic structure
on the foliated manifold $(N^*{\cal F},{\cF}_N)$, and
$r_N^*(|dy\wedge d\eta|^{1/2})\in C^{\infty}(G_{{\cal F}_N},
|H{{\cG}_{{\cal F}_N}}|^{1/2})$ its pull back
via the map $r_N: G_{{\cal F}_N}\rightarrow N^*{\cal F}$.

Recall that the space $S^{m}(G_{{\cal F}_N},
|TG_{{\cal F}_N}|^{1/2})$ is defined to be the space of all
smooth sections $s$
of the vector bundle $|TG_{{\cal F}_N}|^{1/2}$
on $G_{{\cal F}_N}$ homogeneous of degree $m$ such that
$\pi_G({\rm supp}\; s)$ is compact in $G_{\cal F}$.

The half-density principal symbol of $R(k)$ is an element
of $S^{0}(G_{{\cal F}_N}, |TG_{{\cal F}_N}|^{1/2})$
given by the formula
$$
\sigma(R(k))(\gamma,\nu)= \pi^*_Gk(\gamma,\nu)
\otimes r_N^*(|dy\wedge d\eta|^{1/2}), \quad (\gamma,\nu)\in G_{{\cal F}_N}.
$$

The Maslov bundle $M({{\RR}\times G_{{\cal F}_N}},\Pi)$ of
the immersed canonical relation
$({\RR}\times G_{{\cal F}_N},\Pi)$ restricted to $t=0$ is isomorphic to
the Maslov bundle $M(G_{{\cal F}_N})$ of $G_{{\cal F}_N}$, therefore, it
has a canonical constant section, which
extends to a global section $s$ of $M({{\RR}\times G_{{\cal F}_N}},\Pi)$
by requiring it to be constant along each bicharacteristic
$(t,\tau,\nu_1,\nu_2), \nu_1=f_{-t}(\nu_2), t\in {\RR}.$

Using the description of the principal symbol of the operator $W(t)$
given, for instance, in \cite{DG} and the composition theorem of
Fourier integral operators, we get immediately that
the principal symbol of the operator $R(k)W$ is an element of
$S^{0}({\RR}\times G_{{\cal F}_N},
M({{\RR}\times G_{{\cal F}_N}},\Pi)\otimes
|T({\RR}\times G_{{\cal F}_N})|^{1/2})$, whose value at a
point $(t,\gamma,\nu)\in {\RR}\times
G_{{\cal F}_N}$ is given by
\begin{equation*}
\sigma(R(k)W)(t,\gamma,\nu) =
e^{i\int^t_0 \sigma_{\rm sub}(P)(f_{-s}dh_{\gamma}^*(\nu))ds}
s\otimes |dt|^{1/2}\otimes \pi^*_Gk(\gamma,\nu) \otimes
r_{N}^*(|dy\wedge d\eta|^{1/2}).
\end{equation*}
Now let us turn to the composition (\ref{trace}). First, we check
the corresponding cleanness assumption.
\begin{lem}
The assumption on the flow $f_t$ to be clean on ${\cZ}_t$ in
the sense of Definition~\ref{clean} guarantees that the composition
of ${\RR}\times G_{{\cal F}_N}$ with $\Gamma$ is clean.
\end{lem}
\begin{proof}
By definition, the composition of ${\RR}\times G_{{\cal F}_N}$
with $\Gamma $ is clean iff ${\cZ}_t$ is a submanifold
of ${\RR}\times G_{{\cal F}_N}$ and in addition
the fiber product diagram~(\ref{diag}) is clean at any point
$(t,\gamma,\nu)\in {\cZ}$, that is,
the linearized diagram
\begin{equation}
\label{clean:diag}
\CD
T_{p_1(t,\gamma,\nu)}\Gamma @<dp_1<< T_{(t,\gamma,\nu)}{\cZ}\\
@Vd\varphi VV @Vdp_2 VV\\
T_{(t,p(\nu),\nu,-\nu)}(T^*({\RR}\times M\times M))
@<d\Pi << T_{(t,\gamma,\nu)}({\RR}\times G_{{\cal F}_N})
\endCD
\end{equation}
is a fiber product diagram. Since $T_{(t,\gamma,\nu)}{\cZ}$ is always contained
in $T_{\nu}N^*{\cF}\oplus T_{dh^*_{\gamma}(\nu)}{N^*\cF}$, this
is true iff the diagram
\begin{equation}\label{clean:diag1}
\CD
T_{p_1(t,\gamma,\nu)}(\Gamma \bigcap {T^*\RR}\times
N^*{\cF}\times N^*{\cF}\times {T^*\RR}) @<dp_1<< T_{(t,\gamma,\nu)}{\cZ}\\
@Vd\varphi VV @Vdp_2 VV\\
T_{(t,p(\nu),\nu,-\nu)}({T^*\RR}\times
 N^*{\cF}\times N^*{\cF})
@<d\Pi << T_{(t,\gamma,\nu)}({\RR}\times G_{{\cal F}_N})
\endCD
\end{equation}
is a fiber product diagram.

The diagram (\ref{clean:diag1}) has a subdiagram
\begin{equation}\label{leaf:diag}
\CD
L_{p_1(t,\gamma,\nu)}\Gamma @<dp_1<< L_{(t,\gamma,\nu)}{\cZ}\\
@Vd\varphi VV @Vdp_2 VV\\
0 \oplus T_{\nu}{{\cF}_N}\oplus
T_{-\nu}{{\cF}_N}
@<d\Pi << L_{(t,\gamma,\nu)}({{\RR}\times G_{{\cF}_N}})
\endCD
\end{equation}
where $L_{(t,\gamma,\nu)}({{\RR}\times G_{{\cF}_N}})\subset
T_{(t,\gamma,\nu)}({{\RR}\times G_{{\cF}_N}})$ is given by
\begin{multline*}
L_{(t,\gamma,\nu)}({{\RR}\times G_{{\cF}_N}})\\
=\{(U,V_1,V_2,W)\in
T_{t}({\RR})\oplus T_{\nu}{\cF}_N\oplus T_{dh^*_{\gamma}(\nu)}{\cF_N}
 \oplus  H_{(\gamma,\nu)}{{\cG}_{{\cal F}_N}} : U=0, W=0\}\\
\cong T_{\nu}{\cF}_N\oplus T_{dh^*_{\gamma}(\nu)}{\cF_N}
\end{multline*}
$L_{(t,\gamma,\nu)}{\cZ}\subset
T_{(t,\gamma,\nu)}{\cZ}$ by
\begin{multline*}
L_{(t,\gamma,\nu)}{\cZ}
=\{(U,V_1,V_2,W)\in
T_{t}({\RR})\oplus T_{\nu}{\cF}_N\oplus T_{dh^*_{\gamma}(\nu)}{\cF_N}\oplus
H_{(\gamma,\nu)}{{\cG}_{{\cal F}_N}}\\
:U=0,V_1=-df_{-t}(dh^*_{\gamma}(\nu))(V_2),W=0\} \cong T_{\nu}{\cF}_N
\end{multline*}
and $L_{p_1(t,\gamma,\nu)}\Gamma
\subset T_{p_1(t,\gamma,\nu)}\Gamma$ by
\begin{multline*}
L_{p_1(t,\gamma,\nu)}\Gamma \cong \{(U_1,V_1,V_2,U_2)\in
T_{(t,p(\nu))}({T^*\RR})\oplus
T_{\nu}{{\cF}_N}\oplus T_{-\nu}{{\cF}_N} \oplus T_{(t,-p(\nu))}({T^*\RR})\\
: U_1=U_2=0, V_1=-V_2\}
\end{multline*}
which can be easily seen to be a fiber diagram.

Therefore, the diagram \eqref{clean:diag1} is a fiber product diagram
iff the quotient diagram  is a fiber product diagram:
\begin{equation}
\label{trans:diag}
\CD
H_{p_1(t,\gamma,\nu)}\Gamma @<dp_1<< H_{(t,\gamma,\nu)}{\cZ}\\
@Vd\varphi VV @Vdp_2 VV\\
T_{(t,p(\nu))}({T^*\RR})\oplus H_{\nu}{{\cF}_N}\oplus
H_{-\nu}{{\cF}_N}
@<d\Pi << H_{(t,\gamma,\nu)}({{\RR}\times G_{{\cF}_N}})
\endCD
\end{equation}
where
\begin{multline*}
H_{(t,\gamma,\nu)}({{\RR}\times G_{{\cF}_N}})\\
= \{(U,V_1,V_2,W)\in
T_{t}({\RR})\oplus T_{\nu}{\cF}_N\oplus T_{dh^*_{\gamma}(\nu)}{\cF_N}
\oplus  H_{(\gamma,\nu)}{{\cG}_{{\cal F}_N}}: V_1=0,V_2=0\}\\
\cong T_{t}({\RR})\oplus H_{(\gamma,\nu)}{{\cG}_{{\cal F}_N}},
\end{multline*}
\begin{equation*}
H_{(t,\gamma,\nu)}{\cZ}=
T_{(t,\gamma,\nu)}{\cZ}/L_{(t,\gamma,\nu)}{\cZ}
\cong T_{\nu}{Z}_t/T_{\nu}{\cF}_N,
\end{equation*}
\begin{multline*}
H_{p_1(t,\gamma,\nu)}\Gamma
\cong  \{(U_1,V_1,V_2,U_2)\in
T_{(t,p(\nu))}({T^*\RR})\oplus H_{\nu}{{\cF}_N}
\oplus H_{-\nu}{{\cF}_N} \oplus T_{(t,-p(\nu))}({T^*\RR})\\
:U_1=-U_2, V_1=-V_2\}.
\end{multline*}
In its turn, the diagram (\ref{trans:diag}) is a fiber product
diagram iff the flow $f_t$ is clean on ${\cZ}_t$ in the sense of
Definition~\ref{clean}.
\end{proof}

For any connected component ${\cZ}_j$ of ${\cZ}_t$, the excess
of the clean diagram~(\ref{diag}) equals $d_j$, the dimension
of the relative fixed point set $S{\cZ}_j$ in $G_{SN^*{\cF}}$.
By the composition theorem of Fourier integral operators, $\theta_k$ belongs to
$\bigoplus_jI^{\frac{d_j-p}{2}-\frac{1}{4}}(\Lambda_t)$,
where $\Lambda_t=\{(t,\tau)\in T^*{\RR}:
\tau\in {\RR}_{-}\}$,
that proves  the desired representation of $\theta_k$ in the form
(\ref{Ip}) and the existence of the asymptotic expansion~(\ref{expansion})
for $\alpha_j(s,k)$.\footnote{Note that the appearance of the term
$-p/2$ in the exponent is due to the fact
$R(k)W\in I^{-p/2-1/4}(M\times M;\RR\times G_{{\cal F}_N},\Pi)$.}

To obtain the explicit formula for the leading coefficients $\alpha_{j,0}$, we
 compute the principal symbol of $\theta_k$, $\sigma(\theta_k)$, following the
arguments in \cite{Gu-Uribe}.

Fix a connected component ${\cZ}_j$ of ${\cZ}_t$ and $(t,\gamma,\nu)\in
{\cZ}_j$.
The fiber product diagram \eqref{clean:diag} defines
a composition map \cite{DG,H4}
\begin{equation}
* : |T_{p_1(t,\gamma,\nu)}\Gamma|^{1/2}\otimes
|T_{(t,\gamma,\nu)}({\RR}\times G_{{\cal F}_N})|^{1/2}
 \rightarrow |T_{(t,-p(\nu))}(\Lambda_t)|^{1/2}
 \otimes   |T_{(t,\gamma,\nu)}{\cZ}_j|
\end{equation}
and, due to (\ref{trace}) and the composition theorem
for Fourier integral operators, the value of the principal symbol
$\sigma(\theta_k)\in |T\Lambda_t|^{1/2}$ at a point
$(t,\tau)\in \Lambda_t$ is given by
integration over ${\cZ_j}$
of $\sigma(\pi_*\Delta^*)\ast\sigma(W_k)\in
C^{\infty}(\Lambda_t\times \cZ_j, |T\Lambda_t|^{1/2}\otimes |T\cZ_j|)$.

Using the functoriality of the $\ast$ operation on half densities with respect
to reduction (cf. \cite[proof of Lemma 4.6]{Gu-Uribe}),
the  computation of the $\ast$-product $\sigma(\pi_*\Delta^*)\ast\sigma(W_k)$
can be reduced to the computation of a $\ast$-product defined by
the transverse fiber product diagram (\ref{trans:diag}):
\begin{equation*}
\ast_t :|H_{p_1(t,\gamma,\nu)}\Gamma|^{1/2} \otimes
|H_{(t,\gamma,\nu)}({\RR}\times G_{{\cal F}_N})|^{1/2}
 \rightarrow  |T_{(t,-p(\nu))}\Lambda_t|^{1/2}
\otimes |H_{(t,\gamma,\nu)}{\cZ}_j|.
\end{equation*}

More precisely, we apply the result stated in
\cite[proof of Lemma 4.6]{Gu-Uribe} with a symplectic vector space
$\cV=T_{(t,p(\nu),\nu,-\nu,t,-p(\nu))}(T^*{\RR}\times
T^*M\times T^*M\times T^*{\RR})$,
two Lagrangian subspaces in $\cV$:
$\Lambda_1=T_{p_1(t,\gamma,\nu)}\Gamma$ and
$\Lambda_2$, which is the image of $T_{(t,\gamma,\nu)}({\RR}\times
G_{{\cal F}_N})$
in $\cV$ and the reduction given by the coisotropic subspace
$\Gamma=T_{(t,p(\nu),\nu,-\nu,t,-p(\nu))}(T^*{\RR}\times N^*\cF\times
N^*\cF\times T^*{\RR}).$

By this result, the leafwise component of $\sigma(\theta_k)$,
$\sigma_l(\theta_k)\in |L_{(t,\gamma,\nu)}{\cZ}_j|$, is obtained from
the leafwise component of $\sigma(R(k)W)$:
$$
\sigma_l(R(k)W)=
e^{i\int^t_0 \sigma_{\rm sub}(P)(f_{-s}dh_{\gamma}^*(\nu))ds}
\pi^*_Gk
$$
by application of the restriction map $R_{\cZ}$:
$$
\sigma_l(\theta_k)=e^{i\int^t_0
\sigma_{\rm sub}(P)(f_{-s}dh_{\gamma}^*(\nu))ds}
R_{\cZ}\pi^*_Gk,
$$
and the transversal component of $\sigma(\theta_k)$,
$\sigma_t(\theta_k)\in |T_{(t,-p(\nu))}(\Lambda_t)|^{1/2}
\otimes |H_{(t,\gamma,\nu)}{\cZ}_j|$, is equal to
the $\ast_t$-product of the transverse component of $\sigma(R(k)W)$,
$$\sigma_t(R(k)W)=|dt|^{1/2}\otimes r_{N}^*(|dy\wedge d\eta|^{1/2}),$$
and the transverse component of $\sigma(\pi_*\Delta^*)$.
By \cite{DG}, we get
$$
\sigma_t(\theta_k)=|dt|^{1/2}\otimes d\mu_{\cZ_j}
$$
that completes the calculation of the half-density principal symbol
of $\theta_k$ and
implies the formula~\eqref{e:leading}
for the leading coefficients $\alpha_{j,0}$ as in
\cite{DG}.

\section{Maslov indices}
In this section, we define the Maslov factors $\sigma$, corresponding to
$(\gamma,\nu)\in {\cal Z}_t$. For this goal, we will use local
coordinates on the holonomy groupoid $G$ described,
for instance, in \cite{Co79,noncom}.
Choose a pair of compatible foliated charts
near the points $\pi(\nu)$ and $\pi(f_{t}(\nu))$ with the coordinates
$(x,y)$ and $(x',y)$, corresponding to $\gamma\in G$.
Then we have the corresponding coordinates in
$H_{\nu}{{\cal F}_N}=T_{\nu}N^*{\cal F}/T_{\nu}{{\cal F}_N}$
defined as $(\delta y,\delta\eta)$,
and the vertical and horizontal subspaces $V_{\nu}$ and $H_{\nu}$, given
by the equations $\delta y=0$ and $\delta\eta =0$ accordingly.
The linear holonomy map $dH_{(\gamma,\nu)}$
of the foliation $(N^*\cF,\cF_N)$ defines an
isomorphism of the symplectic spaces
$H_{f_t(\nu)}{{\cal F}_N}$ and $H_{\nu}{{\cal F}_N}$,
which preserves the vertical and horizontal subspaces.
Due to this isomorphism, we can obtain a closed curve
$\omega_{(\gamma,\nu)}$ in the Lagrangian Grassmannian ${\cal G}$
of the symplectic space $H_{\nu}{{\cal F}_N}$,
pulling back, via $df_t$, the vertical subspace at $f_t(\nu)$ for $t$
between $0$ and $T$. Denote by $\kappa_{(\gamma,\nu)}$ the
intersection number of $\omega$ with the horizontal subspace $H_{\nu}$:
$$
\kappa_{(\gamma,\nu)}=[\omega_{(\gamma,\nu)}:H_{\nu}].
$$

Let $\chi(t,x,y,\xi,\eta)$ be the generating function of the canonical
transformation $f_t$ in the chosen coordinates. Recall that $\chi$ is
the solution of the Cauchy problem
\begin{equation}
\label{gen}
d_{t}\chi = p(x,y,d_x\chi,d_y\chi), \quad
\chi(0,x,y,\xi,\eta)=x\xi+y\eta.
\end{equation}
By the holonomy invariance of $p$,
it can be easily seen that $\chi(t,x,y,0,\eta)$ is independent of $x$
and $\xi$: $\chi(t,x,y,0,\eta)=\chi(t,y,\eta)$.
Let
\[
R_{(\gamma,\nu)}=\left[
\begin{array}{ccc}
d^2_{yy}\chi & d^2_{y\eta}\chi & -1\\
d^2_{\eta y}\chi & d^2_{\eta\eta}\chi & 0\\
-1& 0& 0
\end{array}
\right].
\]
We define the Maslov factor $\sigma(\gamma,\nu)$ as
$$
\sigma(\gamma,\nu)={\rm sgn}\;R_{(\gamma,\nu)} +
2\kappa_{(\gamma,\nu)}, \quad (\gamma,\nu)
\in {\cal Z}.
$$
It is clear that $\sigma(\gamma,\nu)$ is a locally constant function on
${\cal Z}$.

To handle with Maslov factors in the proof of Theorem~\ref{t:main},
let us write the Schwartz kernel of the
operator $R(k)$ in a foliated coordinate chart on $G$ given by a
pair of compatible coordinate systems as $k(x,x_1,y)\delta(y-y_1)$,
and the Schwartz kernel $W(t)$ microlocally as an
oscillatory integral of the form
$$
\int e^{i\alpha(t,x_1,y_1,x_2,y_2,\xi,\eta)} a(t,x_1,y_1,x_2,y_2,\xi,\eta)
d\xi d\eta,
$$
where
$$
\alpha (t,x_1,y_1,x_2,y_2,\xi,\eta)=\chi(t,x_1,y_1,\xi,\eta) -
x_2\xi -y_2\eta
$$
and $\chi(t,x_1,y_1,\xi,\eta)$ is the generating function of the
canonical transformation $f_t$ given by (\ref{gen}). Then the Schwartz
kernel of the operator $R(k)W$ is given by the formula
$$
\int e^{i\alpha(t,x_1,y_1,x_2,y_2,\xi,\eta)}k(x,x_1,y_1)
a(t,x_1,y_1,x_2,y_2,\xi,\eta)
dx_1 d\xi d\eta,
$$
from where one can easily derive the desired assertion, following
the arguments of \cite{DG}.

\pagebreak


\begin{thebibliography}{10}

\bibitem{CdV73}
Y.~Colin de~Verdi\`ere.
\newblock Spectre du {Laplacien} et longeurs des g\'eod\'esiques periodiques.
\newblock {\em Comp. Math.}, 27:159--184, 1973.

\bibitem{Chaz}
J.~Chazarain.
\newblock Formule de {Poisson} pour le vari\'{e}t\'{e}s riemanniennes.
\newblock {\em Invent. Math.}, 24:65--82, 1974.

\bibitem{DG}
J.J. Duistermaat and V.~Guillemin.
\newblock The spectrum of positive elliptic operators and periodic
  bicharacteristics.
\newblock {\em Invent.Math.}, 29:39--79, 1975.

\bibitem{noncom}
Yu.~A. Kordyukov.
\newblock Noncommutative spectral geometry of {Riemannian} foliations.
\newblock {\em Manuscripta Math.}, 94:45--73, 1997.

\bibitem{Connes:zeta}
A.~Connes.
\newblock Trace formula in noncommutative geometry and the zeros of the
  {Riemann} zeta function.
\newblock {\em Selecta Math. (N.S.)}, 5:29--106, 1999.

\bibitem{Golse:Leichtnam}
F.~Golse and E.~Leichtnam.
\newblock Applications of {Connes'} geodesic flow to trace formulas in
  noncommutative geometry.
\newblock {\em J. Funct. Anal.}, 160:408--436, 1998.

\bibitem{Re}
B.~L. Reinhart.
\newblock {\em Differential Geometry of Foliations}.
\newblock Springer, Berlin Heidelberg New York, 1983.

\bibitem{Co79}
A.~Connes.
\newblock Sur la th\'eorie non-commutative de l'int\'egration.
\newblock In {\em Alg\`ebres d'op\'erateurs}, volume 725 of {\em Lecture Notes
  Math.}, pages 19--143. Springer, Berlin Heidelberg New York, 1979.

\bibitem{GS79}
V.~Guillemin and S.~Sternberg.
\newblock Some problems in integral geometry and some related problems in
  microlocal analysis.
\newblock {\em Amer. J. Math.}, 101:915--959, 1979.

\bibitem{geo}
Yu.~A. Kordyukov.
\newblock The transversal wave equation and the noncommutative geodesic flow in
  {R}iemannian foliations.
\newblock Preprint ETH Z\"urich, 1997; dg-ga/9703015.

\bibitem{H4}
L.~H\"{o}rmander.
\newblock {\em The analysis of linear partial differential operators {I}{V}}.
\newblock Springer, Berlin Heidelberg New York, 1986.

\bibitem{geom:asymp}
V.~Guillemin and S.~Sternberg.
\newblock {\em Geometric Asymptotics}.
\newblock American Mathematical Society, Providence, R. I., 1977.

\bibitem{Gu-Uribe}
V.~Guillemin and A.~Uribe.
\newblock Circular symmetry and the trace formula.
\newblock {\em Invent. Math.}, 96:385--423, 1989.

\end{thebibliography}
\end{document}